\documentclass[a4paper,reqno]{amsart}

\usepackage{amssymb,upref}
\usepackage[mathcal]{euscript}
\usepackage[cmtip,all]{xy}
\usepackage{color}

\newcommand{\bydef}{:=}

\newcommand{\id}{\mathrm{id}}

\newcommand{\diag}{\mathrm{diag}}

\newcommand{\tr}{\mathrm{tr}}

\newcommand{\cA}{\mathcal{A}}
\newcommand{\cB}{\mathcal{B}}

\newcommand{\cD}{\mathcal{D}}

\newcommand{\cH}{\mathcal{H}}

\newcommand{\cK}{\mathcal{K}}
\newcommand{\cL}{\mathcal{L}}

\newcommand{\cQ}{\mathcal{Q}}
\newcommand{\cR}{\mathcal{R}}

\newcommand{\cV}{\mathcal{V}}

\newcommand{\frg}{{\mathfrak g}}

\newcommand{\frd}{{\mathfrak d}}
\newcommand{\frf}{{\mathfrak f}}

\newcommand{\ZZ}{\mathbb{Z}}

\newcommand{\CC}{\mathbb{C}}

\newcommand{\OO}{\mathbb{O}}
\newcommand{\FF}{\mathbb{F}}

\newcommand{\Albert}{\mathbb{A}}

\DeclareMathOperator{\Hom}{\mathrm{Hom}}
\DeclareMathOperator{\End}{\mathrm{End}}

\DeclareMathOperator{\Aut}{\mathrm{Aut}}
\DeclareMathOperator{\Antiaut}{\mathrm{Antiaut}}
\DeclareMathOperator{\Diag}{\mathrm{Diag}}
\DeclareMathOperator{\AAut}{\mathbf{Aut}}

\newcommand{\ad}{\mathrm{ad}}
\newcommand{\Ad}{\mathrm{Ad}}

\newcommand{\frsl}{{\mathfrak{sl}}}
\newcommand{\frsp}{{\mathfrak{sp}}}
\newcommand{\frso}{{\mathfrak{so}}}
\newcommand{\frpsl}{{\mathfrak{psl}}}

\newcommand{\GL}{\mathrm{GL}}

\newcommand{\PGL}{\mathrm{PGL}}

\newtheorem{theorem}{Theorem}

\theoremstyle{definition}

\newtheorem{example}[theorem]{Example}
\newtheorem{remark}[theorem]{Remark}

\numberwithin{theorem}{section}
\numberwithin{equation}{section}

\begin{document}

\title{An overview  of  
fine gradings on simple Lie algebras}

\author[C.~Draper]{Cristina Draper${}^{*}$}
\address{Departamento de Matem\'atica Aplicada, Escuela de las Ingenier\'{\i}as, 
Universidad de M\'alaga, Ampliaci\'on Campus de Teatinos, 29071 M\'alaga, Spain}
\email{cdf@uma.es}
\thanks{${}^{*}$ Supported by the Spanish Ministerio de Econom\'{\i}a y Competitividad---Fondo 
Europeo de Desarrollo Regional (FEDER)  MTM2013-41208-P, and by the Junta de Andaluc\'{\i}a 
grants FQM-336 and FQM-7156, with FEDER funds}

\author[A.~Elduque]{Alberto Elduque${}^{**}$}
\address{Departamento de Matem\'{a}ticas
 e Instituto Universitario de Matem\'aticas y Aplicaciones,
 Universidad de Zaragoza, 50009 Zaragoza, Spain}
\email{elduque@unizar.es}
\thanks{${}^{**}$ Supported by the Spanish Ministerio de Econom\'{\i}a y Competitividad---Fondo 
Europeo de Desarrollo Regional (FEDER) MTM2013-45588-C3-2-P, and by the Diputaci\'on General de 
Arag\'on---Fondo Social Europeo (Grupo de Investigaci\'on de \'Algebra)}

\subjclass[2010]{Primary 17B70; Secondary 17B40, 17B20, 17B60}

\keywords{Simple Lie algebra; grading; fine; quasitorus.}

\date{}

\begin{abstract}
This paper presents a survey of the results and ideas behind the classification of the fine gradings, up to equivalence, on the simple finite dimensional Lie algebras over an algebraically closed field of characteristic zero.
\end{abstract}

\maketitle



\section{Introduction}\label{se:intro}

These notes constitute an expanded version of some parts of the mini course delivered by the second author at the Conference ``Advances in Group Theory and Applications AGTA-2015''.

Gradings by abelian groups have played a key role in the study of Lie algebras and superalgebras, starting with the root space decomposition of the semisimple Lie algebras over the complex field, which is an essential ingredient in the Killing-Cartan classification of these algebras. Gradings by a cyclic group appear in the connection between Jordan algebras and Lie algebras through the Tits-Kantor-Koecher construction, and in the theory of Kac-Moody Lie algebras. Gradings by the integers or the integers modulo $2$ are ubiquitous in Geometry.

In 1989, Patera and Zassenhaus \cite{PZ89} began a systematic study of gradings by abelian groups on Lie algebras. They raised the problem of classifying the fine gradings, up to equivalence, on the simple Lie algebras over the complex numbers. This problem has been settled now thanks to the work of many colleagues.

A list of the corresponding maximal quasitori for the classical simple Lie algebras over $\CC$, with the exception of $D_4$ later considered in \cite{DMV10}, was obtained in \cite{HPP98}, while a full classification, including $D_4$, was given in \cite{Eld10}, relying on the previous work of several other authors.

As for the exceptional simple Lie algebras, fine gradings on $G_2$ were classified independently in \cite{DM06} and \cite{BT09}, based on the classification of gradings on the octonions in \cite{Eld98}. Fine gradings on $F_4$ were classified in \cite{DM09} (see also \cite{Cris_NonCompu}), where these were used to classify the fine gradings on the Albert algebra. The process can be reversed, first classifying the fine gradings on the Albert algebra and then using these to classify the fine gradings on $F_4$, in a way valid also in prime characteristic \cite{EK12}. For $E_6$ the classification of fine gradings was obtained in \cite{DV_E6}. The classification for $E_7$ and $E_8$ can be extracted from the recent work of Jun Yu \cite{Yu12,Yu14}.

However, someone looking for the first time at the problem of the classification of fine gradings on the finite dimensional simple Lie algebras over an algebraically closed field of characteristic zero finds it  difficult to get a unified list of the fine gradings and a neat idea of what the fine gradings look like, without going through many technical details scattered through different places.

The goal of this survey is to provide a description of the classification of fine gradings on the finite dimensional simple Lie algebras over an algebraically closed field $\FF$ of characteristic zero. No proofs will be given but the main ideas used in the classification will be explained.

The reader is referred to the monograph \cite{EK13} and the references  therein for most of the missing details.

\smallskip

The paper is organized as follows. Section \ref{se:background} will be devoted to give the basic definitions on gradings needed in the sequel, as well as the relationship of fine gradings with maximal quasitori of the automorphism group. Then the fine gradings on the simple Lie algebras of type $A$ will be treated in Section \ref{se:A} and on the orthogonal and symplectic Lie algebras in Section \ref{se:BCD}. Fine gradings on the exceptional simple Lie algebras will be quickly reviewed in Section \ref{se:exceptional}. The reader may consult \cite{DE_Note}. Finally, Section \ref{se:modular} will show how most of the results in characteristic zero remain valid in the modular case, and will highlight what remains to be done.

\bigskip


\section{Basic definitions}\label{se:background}

Let $\FF$ be an arbitrary ground field. Vector spaces and algebras will be defined over $\FF$. Unadorned tensor products will be assumed to be over $\FF$.

Given an abelian group $G$ and a nonasociative (i.e., not necessarily associative) algebra $\cA$, a \emph{grading on $\cA$ by $G$} (or \emph{$G$-grading}) is a decomposition into a direct sum of subspaces
\begin{equation}\label{eq:grading}
\Gamma:\cA=\bigoplus_{g\in G}\cA_g
\end{equation}
such that $\cA_{g_1}\cA_{g_2}\subseteq \cA_{g_1g_2}$ for any $g_1,g_2\in G$. For each $g\in G$, $\cA_g$ is the \emph{homogeneous component of degree $g$}, its elements are the \emph{homogeneous elements of degree $g$}.

Assume that
\[
\Gamma_1:\cA=\bigoplus_{g\in G}\cA_g\quad\text{and}\quad\Gamma_2:\cA=\bigoplus_{h\in H}\cA'_h
\]
are two gradings on $\cA$. Then:
\begin{itemize}
\item 
$\Gamma_1$ and $\Gamma_2$ are said to be \emph{equivalent} if there is an automorphism $\varphi\in\Aut(\cA)$, called a \emph{graded equivalence}, such that for any $g\in G$ there is an $h\in H$ with $\varphi(\cA_g)=\cA'_h$.

\item
$\Gamma_1$ is said to be a \emph{refinement} of $\Gamma_2$ if for any $g\in G$ there is an $h\in H$ such that $\cA_g\subseteq \cA'_h$. Then $\Gamma_2$ is said to be a \emph{coarsening} of $\Gamma_1$. If one of these containments is strict, the refinement is said to be \emph{proper}. 
\end{itemize}

The grading $\Gamma$ is said to be \emph{fine} if it admits no proper refinement \cite[Definition 2]{PZ89}. Any grading is a coarsening of a fine grading.

\begin{example}[Cartan grading]\label{ex:Cartan}
Let $\cL$ be a finite dimensional semisimple Lie algebra of rank $r$ over an algebraically closed field of characteristic zero. Let $\cH$ be a Cartan subalgebra of $\cL$ with root system $\Phi$. The root space decomposition
\[
\cL=\cH\oplus\left(\bigoplus_{\alpha\in\Phi}\cL_\alpha\right)
\]
is an example of a fine grading by the group $\ZZ\Phi$ (isomorphic to $\ZZ^r$). Here $\cL_0=\cH$.
\end{example}

\begin{example}[Pauli grading]\label{ex:Pauli}
Assume that $n\geq 2$ and $\FF$ contains a primitive $n$th root of unity $\epsilon$, and let $\cA=M_n(\FF)$ be the associative algebra of $n\times n$ matrices over $\FF$. Consider the matrices
\[
x=\begin{pmatrix}
1&0&0&\ldots&0\\
0&\epsilon&0&\ldots&0\\
0&0&\epsilon^2&\ldots&0\\
\vdots&\vdots&\vdots&\ddots&\vdots\\
0&0&0&\ldots&\epsilon^{n-1}\end{pmatrix},
\qquad
y=\begin{pmatrix}
0&1&0&\ldots&0\\
0&0&1&\ldots&0\\
\vdots&\vdots&\vdots&\ddots&\vdots\\
0&0&0&\ldots&1\\
1&0&0&\ldots&0
\end{pmatrix}.
\]
They satisfy $x^n=y^n=1$, $yx=\epsilon xy$. Then the decomposition
\[
\cA=\bigoplus_{(\bar\imath,\bar\jmath)\in\ZZ_n^2}\cA_{(\bar\imath,\bar\jmath)},\quad\text{with $\cA_{(\bar\imath,\bar\jmath)}=\FF x^iy^j$},
\]
is a fine grading on $\cA$. Moreover, $\cA$ becomes a \emph{graded division algebra}, that is, any nonzero homogeneous element is invertible.

The Pauli grading induces a fine $\ZZ_n^2$-grading in the special linear Lie algebra $\frsl_n(\FF)=\{x\in\cA: \tr(x)=0\}$.
\end{example}

The Cartan grading and the Pauli grading on $\frsl_n(\FF)$ are quite different in nature. For the Cartan grading any homogeneous element of degree $\ne 0$ is $\ad$-nilpotent, while in the Pauli grading any nonzero homogeneous element is $\ad$-semisimple.

It should be noted that the Cartan grading on $\frsl_n(\FF)$ is the restriction of the grading on $M_n(\FF)$, also called Cartan grading, by $\ZZ^{n-1}=\ZZ\epsilon_1\oplus\cdots\oplus\ZZ\epsilon_{n-1}$, where $\epsilon_i=(0,\ldots,1,\ldots,0)$ ($1$ in the $i$th position), such that the degree of $E_{ij}$ is $\epsilon_{i-1}-\epsilon_{j-1}$, where $E_{ij}$ is the matrix with $1$ in the $(i,j)$ slot and $0$'s elsewhere, and where $\epsilon_0=0$. We may think of $M_n(\FF)$ as the algebra of endomorphisms $\End_\FF(\cV)$, where $\cV$ is an $n$-dimensional vector space with a basis $\{v_1,\ldots,v_n\}$. Then $\cV$ is $\ZZ^{n-1}$-graded: $\cV=\bigoplus_{i=0}^{n-1} \cV_{\epsilon_i}$ (simply a decomposition as a direct sum of vector subspaces), with $\cV_{\epsilon_i}=\FF v_{i+1}$ for $i=0,\ldots,n-1$. The Cartan grading is the grading induced on $\End_\FF(\cV)$:
\[
\End_\FF(\cV)_\epsilon\bydef \{ f\in\End_\FF(\cV): f(\cV_\delta)\subseteq \cV_{\epsilon+\delta}\ \forall\delta\in\ZZ^{n-1}\}.
\]

\begin{remark}\label{re:Qm}
The case of $n=2$ in Example \ref{ex:Pauli} will appear quite often in what follows. For further reference consider the matrices
\begin{equation}\label{eq:qs}
q_1=\begin{pmatrix} 1&0\\ 0&-1\end{pmatrix},\quad 
q_2=\begin{pmatrix} 0&1\\ 1&0\end{pmatrix},\quad
q_3=q_1q_2=\begin{pmatrix} 0&1\\ -1&0\end{pmatrix}.
\end{equation}
Denote by $\cQ$  the algebra $M_2(\FF)$ with the $\ZZ_2^2$-grading (the Pauli grading) where $q_1$ is homogeneous of degree  $(\bar 1,\bar 0)$ and  $q_2$ is homogeneous of degree $(\bar 0,\bar 1)$. This is a fine grading and the transpose involution preserves the homogeneous components.

Any tensor power $\cQ^{\otimes m}$ is endowed with the naturally induced grading by $\ZZ_2^{2m}\,\bigl(\simeq (\ZZ_2^2)^m\bigr)$. This is a division grading, and the involution $\tau$  which acts as the transpose on each factor  is an orthogonal involution that preserves each homogeneous component. (If $m=0$, $\tau=\id$.)
\end{remark}

\smallskip

Given a grading $\Gamma$ as in \eqref{eq:grading} and a character $\chi$ of $G$ (i.e., a group homomorphism $G\rightarrow \FF^\times$), the map $\varphi_\chi:\cA\rightarrow \cA$, such that $\varphi_\chi(x)=\chi(g)x$ for any $g\in G$ and $x\in\cA_g$, is an automorphism of $\cA$ that acts diagonally on $\cA$, as it acts as a scalar on each homogeneous component.

Hence each $\varphi_\chi$ belongs to the \emph{diagonal group} of $\Gamma$, defined as follows:
\[
\Diag(\Gamma)\bydef\{ \varphi\in\Aut(\cA): \forall g\in G,\ \exists \alpha_g\in\FF^\times\ \text{such that}\ \varphi\vert_{\cA_g}=\alpha_g\id\}.
\]
This is a subgroup of $\Aut(\cA)$ (closed in the Zariski topology).

If $\FF$ is algebraically closed of characteristic zero, characters separate elements of $G$, and hence the homogeneous components are just the common eigenspaces for the action of the subgroup $\{\varphi_\chi :\chi\in\hat G\}$, where $\hat G$ denotes the group of characters of $G$. Conversely, assume that $\cA$ has finite dimension and let $K$ be an abelian subgroup of $\Aut(\cA)$ whose elements act diagonally on $\cA$. The common eigenspaces of the action of the elements in $K$ give a grading on $\cA$ by the group of characters of the Zariski closure of $K$. In particular, gradings by $G$ on $\cA$ correspond bijectively to homomorphisms (as algebraic groups) $\hat G\rightarrow \Aut(\cA)$.

The next result \cite[Theorem 2]{PZ89} follows easily:

\begin{theorem}\label{th:gradings_MADs}
Let $\cA$ be a finite dimensional algebra over an algebraically closed field $\FF$ of characteristic zero. Then a grading $\Gamma$ as in \eqref{eq:grading} is fine if and only if $\Diag(\Gamma)$ is a maximal abelian diagonalizable subgroup (i.e., a maximal quasitorus) of $\Aut(\cA)$.

Moreover, two fine gradings on $\cA$ are equivalent if and only if the corresponding diagonal groups are conjugate in $\Aut(\cA)$, so there is a bijection
\begin{equation}\label{eq:fineMAD}
\begin{split}
\left\{\begin{matrix} \text{Equivalence classes}\\
  \text{of fine gradings on $\cA$} \end{matrix}\right\}\quad 
  &\longleftrightarrow\quad
  \left\{\begin{matrix} \text{Conjugacy classes of}\\
  \text{maximal quasitori of $\Aut(\cA)$}
  \end{matrix}\right\}\\
  [\Gamma]\qquad\quad&\longrightarrow\qquad\quad[\Diag(\Gamma)].
\end{split}
\end{equation}
\end{theorem}

As a direct consequence, if $\cA$ and $\cB$ are finite dimensional algebras over an algebraically closed field $\FF$ of characteristic $0$ such that $\Aut(\cA)$ and $\Aut(\cB)$ are isomorphic algebraic groups, then the problems of classifying fine gradings on $\cA$ and on $\cB$ are equivalent.

In Example \ref{ex:Cartan}, the maximal quasitorus attached to the Cartan grading is the maximal torus consisting of those automorphisms of $\cL$ that fix all the elements of the Cartan subalgebra $\cH$. This is isomorphic to $(\FF^\times)^r$. In Example \ref{ex:Pauli} the corresponding maximal quasitorus is the subgroup generated by $\Ad_x$ and $\Ad_y$ (where $\Ad_p(q)\bydef pqp^{-1}$) in $\Aut(M_n(\FF))\simeq \PGL_n(\FF)$.

\smallskip

A grading $\Gamma$ as in \eqref{eq:grading} may be realized by different groups. Think, for instance, of the trivial grading $\cA=\cA_e$, which is a grading by any abelian group. However, there is always a natural grading group: the group of characters of $\Diag(\Gamma)$ (i.e., of homomorphisms of algebraic groups $\Diag(\Gamma)\rightarrow \FF^\times$). This is  called  the \emph{universal group}, or universal grading group. (See \cite[\S 1.4]{EK13}.)

\smallskip

The definition of grading on an algebra $\cA$ admits natural generalizations. For instance, let $\varphi$ be an involution of $\cA$, that is, an involutive antiautomorphism of $\cA$. Then a grading on $(\cA,\varphi)$ is a grading on the algebra $\cA$ as in \eqref{eq:grading} such that $\varphi(\cA_g)=\cA_g$ for any $g\in G$. If $\Aut(\cA,\varphi)$ denotes the group of automorphisms of $\cA$ that commute with $\varphi$, then the bijection in \eqref{eq:fineMAD} becomes a bijection:
\begin{equation}\label{eq:fineMAD_Aphi}
\begin{split}
\left\{\begin{matrix} \text{Equivalence classes of}\\
  \text{fine gradings on $(\cA,\varphi)$} \end{matrix}\right\}\quad 
  &\longleftrightarrow\quad
  \left\{\begin{matrix} \text{Conjugacy classes of}\\
  \text{maximal quasitori of $\Aut(\cA,\varphi)$}
  \end{matrix}\right\}\\
  [\Gamma]\qquad\quad&\longrightarrow\qquad\quad[\Diag(\Gamma)].
\end{split}
\end{equation}
This works too for an antiautomorphism $\varphi$, not necessarily involutive.

\bigskip


\section{Fine gradings on the special linear Lie algebras}\label{se:A}

In this and the next two sections, the ground field $\FF$ will be assumed to be algebraically closed of characteristic zero.

The group $\Aut(\frsl_n(\FF))$ ($n\geq 2$) is determined as follows \cite{Jac79}:

\begin{itemize}
\item
Any automorphism of $\frsl_2(\FF)$ is the restriction of an automorphism of $M_2(\FF)$, so we have an isomorphism $\Aut(\frsl_2(\FF))\simeq \Aut(M_2(\FF))$.

\item 
Any automorphism of $\frsl_n(\FF)$, $n\geq 3$, is the restriction of either an automorphism of $M_n(\FF)$ or the negative of an antiautomorphism of $M_n(\FF)$, so
\[
\Aut(\frsl_n(\FF))\simeq \Aut(M_n(\FF))\cup \bigl(-\Antiaut(M_n(\FF))\bigr).
\]
\end{itemize}
We will identify $\Aut(\frsl_n(\FF))$ with the corresponding group in $\GL\bigl(M_n(\FF)\bigr)$.

Hence, given any maximal quasitorus $M$ in $\Aut(\frsl_n(\FF))$, either:
\begin{itemize}
\item $M\subseteq \Aut(M_n(\FF))$, so that the corresponding fine grading on $\frsl_n(\FF)$ is the restriction of a fine grading on $M_n(\FF)$ (this is always the case for $n=2$); or

\item
there is an antiautomorphism $\varphi$ of $M_n(\FF)$ and a quasitorus $M'$ of $\Aut(M_n(\FF))$ such that $M=M'\cup M'(-\varphi)$. In particular, $\varphi^2\in M'$. In this case, by maximality, $M'$ is a maximal quasitorus of $\Aut\bigl(M_n(\FF),\varphi\bigr)$.
\end{itemize}

In the first case, our task is to find the fine gradings on $M_n(\FF)$, and this is relatively simple. The classical Wedderburn theory tells us that any finite dimensional central simple associative algebra is, up to isomorphism, the algebra of endomorphisms of a finite dimensional right vector space over a central division algebra. The same arguments (see \cite{BSZ01}) imply  that given any grading on $\cR=M_n(\FF)$, $\cR$ is, up to graded isomorphism, the algebra of endomorphisms of a graded right free-module of finite rank over a graded central division algebra: $\cR\simeq \End_\cD(\cV)$.

The graded central division algebras are easily shown to be tensor products of matrix algebras with Pauli gradings, and their degrees can be taken to be powers of prime numbers (see, for instance, \cite[Proposition 2.1]{Eld10}): 
\[
\cD\simeq M_{n_1}(\FF)\otimes\cdots \otimes M_{n_r}(\FF),
\] 
where $n_1,\ldots,n_r$ are powers of prime numbers, with each slot endowed with the Pauli grading as in Example \ref{ex:Pauli}.

On the other hand, if $M=M'\cup M'(-\varphi)$ for an antiautomorphism $\varphi$, we get some freedom as we may change $\varphi$ by $\psi\varphi$ for any $\psi\in M'$. The antiautomorphism $\varphi$ induces an involution preserving the grading in the graded division algebra determined as above for $M'$. But note that if $x$ and $y$ are homogeneous elements with $yx=\epsilon xy$ with $\epsilon^m=1$, and if $\tau$ is an involution that preserves the one-dimensional homogeneous components, then from $\tau(xy)=\tau(y)\tau(x)$, we also get $xy=\epsilon yx$, so that $yx=\epsilon^2 yx$ and $\epsilon^2=1$. This shows that, in this case, $n_1=\cdots=n_r=2$, so our graded division algebra $\cD$ must be isomorphic to $\cQ^{\otimes m}$ (see Remark \ref{re:Qm}) for some $m\geq 0$. For $m=0$, $\cD$ is the ground field $\FF$. Moreover, the involution of $\cD$ can always be taken to be the involution $\tau$ in Remark \ref{re:Qm}.

It turns out that, identifying $\cR$ with $\End_\cD(\cV)$, the antiautomorphism $\varphi$ becomes the `adjoint' relative to a nondegenerate sesquilinear form $B:\cV\times\cV\rightarrow \cD$. That is, $B$ is $\FF$-bilinear, $\cD$-linear in the second component, $B(vd,w)=\tau(d)B(v,w)$, and $B(xv,w)=B(v,\varphi(x)w)$ for any $d\in\cD=\cQ^{\otimes m}$, $x\in \cR$, and $v,w\in \cV$.

By adjusting $\varphi$, using the freedom explained above, and $B$, we may find a homogeneous $\cD$-basis $\{v_1,\ldots,v_r,\ldots,v_{r+2s}\}$  of $\cV$ such that the coordinate matrix of $B$ is of the following block-diagonal form
\begin{equation}\label{eq:matrizdeB} 
M_B=  \scriptstyle{
\begin{pmatrix}
d_1&&&&&&&\\
&\ddots&&&&&&\\
&&d_r&&&&&\\
&&&0&1&&&\\
&&&1&0&&&\\
&&&&&\ddots&&\\
&&&&&&0&1\\
&&&&&&1&0
\end{pmatrix}}
 \end{equation} 
where $r\geq 0$, $d_1,\ldots,d_r$ are nonzero homogeneous elements in $\cD=\cQ^{\otimes m}$, and $\deg(v_i)=g_i$, $i=1,\ldots,r$, with
\begin{equation}\label{eq:gis}
g_i^2=\deg(d_i)\quad \text{for $i=1,\ldots,r$},\qquad g_{r+2i-1}g_{r+2i}=e\quad \text{for $i=1,\ldots,s$.}
\end{equation}
Identifying $\cR\simeq\End_\cD(\cV)$ with $M_{r+2s}(\cD)$, the antiautomorphism $\varphi$ acts on $x=(x_{ij})$ as $\varphi(x)=M_B^{-1}(\tau(x_{ji}))M_B$. Note that the adjoint $\varphi$ needs not be involutive.
 This happens only if $M_B^{-1}M_B^{t\tau}$ is in the center of $M_{r+2s}(\cD)$, which consists of the scalar multiples of the identity matrix (for instance, this is the case for $\cD=\mathbb F$). 

Let $\widetilde G$ be the abelian group generated by a subgroup $T$ isomorphic to $\ZZ_2^{2m}$ (the grading group of $\cD=\cQ^{\otimes m}$) and elements $g_1,\ldots,g_{r+2s}$, subject only to the relations in \eqref{eq:gis}. The universal group of the $\varphi$-grading on $\cR$ with maximal quasitorus $M'$ is isomorphic to the subgroup $\bar G$ of $\widetilde G$ generated by $T$ and the elements $g_ig_j^{-1}$, $1\leq i,j\leq r+2s$. The free rank of $\bar G$ is exactly $s$ and $\bar{G}$ is the cartesian product of a $2$-group that contains $T$ and a free subgroup.

The automorphism $-\varphi$ of $\frsl_n(\FF)$ refines this fine grading on $(\cR,\varphi)$ to a fine grading on $\frsl_n(\FF)$. The universal group of this fine grading (corresponding to the maximal quasitorus $M$) is a group $G$ containing an element $h$ of order $2$ such that $G/\langle h\rangle$ is isomorphic to $\bar G$, since $(-\varphi)^2\in M'$. Recall that the character group $\hat G$ is isomorphic to $M$. The characters of $G$ which are trivial on $h$ correspond to the elements in $M'$ and can be identified with the characters of $\bar G$, while those characters $\chi$ of $G$ with $\chi(h)=-1$ correspond to the elements in $M'(-\varphi)$. For details see \cite{Eld10} or \cite[Chapter 3]{EK13}. If $\varphi$ has order two, then $G$ is isomorphic to $\bar G\times\ZZ_2$.

In the situation above, attach to the maximal quasitorus $M$ of $\Aut(\frsl_n(\FF))$ the sequence  $(m,s;d_1,\ldots,d_r)$ and denote by   $\Gamma_{(m,s;\,d_1,\ldots,d_r)}$ the fine grading above whose diagonal group is $M$. 
 Summing up, for each $m\ge0$ such that $2^m$ divides $n$ and each $s\ge0$ with $s\le n2^{-m-1}$, we take the $\mathbb Z_2^{2m}$-graded division algebra $\cD=\cQ^{\otimes m}$ with the involution $\tau$ acting as the transpose $t$ on each factor, and choose a homogeneous element  $d_i\in\cD$ for each $i\le r=n2^{-m}-2s$. Then, we consider  the  right free-module $\cV$ over $\cD$ with basis $\{v_1,\dots,v_{r+2s}\}$ which is $\widetilde G$-graded by assigning   $\deg(v_i)=g_i$ and imposing  $\cV_g\cD_h\subset\cV_{gh}$ for any $g,h\in \widetilde G$. Thus $\End_\cD(\cV)$ is also $\widetilde G$-graded, where $x\in \End_\cD(\cV)$ has degree $g$ if $xV_h\subset V_{gh}$ for all $h\in \widetilde G$. If we consider the sesquilinear form $B\colon\cV\times\cV\rightarrow \cD$ with coordinate matrix given by   $M_B$ in \eqref{eq:matrizdeB}, we observe that $B$ satisfies $B(\cV_g,\cV_h)=0$ whenever $gh\ne e$, and this implies that the adjoint $\varphi\colon \End_\cD(\cV)\to\End_\cD(\cV)$ relative to $B$ is compatible with the $\widetilde G$-grading on $\End_\cD(\cV)$, whose universal grading group is in fact $\bar G$. The grading $\Gamma_{(m,s;\,d_1,\ldots,d_r)}$ is then the grading on $\frsl_n(\FF)$  considered by restricting the grading on $\End_\cD(\cV)\simeq M_{r+2s}(\cD)\simeq M_n(\FF)$ and refining it with the antiautomorphism $-\varphi$.

Any outer fine grading on $\frsl_n(\FF)$ appears in this way, but not conversely. The grading $\Gamma_{(m,s;\,d_1,\ldots,d_r)}$ is not fine for $s=0$, $r=2$ and $\FF d_1=\FF d_2$. These constitute the only exceptions \cite[Theorem 3.30]{EK13}.

The point is now   to distinguish whether two of these gradings are equivalent.
The equivalence classes of fine gradings on $\frsl_n(\FF)$ are determined as follows:

\begin{theorem}\label{th:A}
\null\quad\null
\begin{enumerate}
\item[(1)] Up to equivalence, the only fine gradings on $\frsl_2(\FF)$ are the Cartan grading (by $\ZZ$) and the Pauli grading (by $\ZZ_2^2$).

\item[(2)] If $n\geq 3$, any fine grading of $\frsl_n(\FF)$ is equivalent to a grading of one and only one of the following types:
\begin{enumerate}
\item[(2.a)] \textup{INNER GRADINGS:} The restriction of a fine grading on $\cR=M_n(\FF)$ obtained by decomposing $n$ as a product $n=mp_1^{s_1}\cdots p_r^{s_r}$, $m\geq 1$, $r\geq 0$, with $p_1,\ldots,p_r$ prime numbers, $s_1,\ldots,s_r\geq 1$; and by considering the isomorphism
\[
M_n(\FF)\simeq M_m(\FF)\otimes M_{p_1^{s_1}}(\FF)\otimes\cdots\otimes M_{p_r^{s_r}}(\FF)
\]
(Kronecker product). The grading on $\cR$ is then given by combining the Cartan grading on $M_m(\FF)$ and the Pauli gradings on each $M_{p_i^{s_i}}(\FF)$. Moreover, if $p_1^{s_1}=\cdots=p_r^{s_r}=2$, we require $m\geq 3$.

\item[(2.b)] \textup{OUTER GRADINGS:} The fine gradings {$\Gamma_{(m,s;\,d_1,\ldots,d_r)}$} with $n=2^m(r+2s)$ and nonzero homogeneous elements $d_i\in\cQ^{\otimes m}$, except for $s=0$, $r=2$, $\FF d_1=\FF d_2$. 

Two such gradings $\Gamma_{(m,s;\,d_1,\ldots,d_r)}$ and $\Gamma_{(m',s';\,d'_1,\ldots,d'_{r'})}$ are equivalent if and only if $m=m'$, $s=s'$, $r=r'$, and there  {exist} a graded equivalence $\Phi:\cQ^{\otimes m}\rightarrow \cQ^{\otimes m}$, a permutation $\sigma\in S_r$ and a nonzero homogeneous element $z\in\cQ^{\otimes m}$ such that $\Phi(zd_i)\in \FF d'_{\sigma(i)}$ for all $i=1,\ldots,r$. 
\end{enumerate}
\end{enumerate}
\end{theorem}

The details can be found in \cite{Eld10} or \cite[Chapter 3]{EK13}. The last requirement in (2.a) is due to the fact that if   $p_1^{s_1}=\cdots=p_r^{s_r}=2$ and $m=1$ or $m=2$, the grading can be refined to a grading in (2.b).

\begin{example}[Fine gradings on $\frsl_4(\FF)$]
\null\quad\null
\begin{itemize}
\item[\textbf{Inner:}] Up to equivalence there are the following possibilities:
\begin{itemize}
\item The Cartan grading by $\ZZ^3$.
\item The Pauli grading by $\ZZ_4^2$.
\end{itemize}

\item[\textbf{Outer:}] In this case, the possible sequences up to equivalence are the following:
\begin{itemize}
\item $(0,2;\emptyset)$, which gives a fine grading by $\ZZ^2\times\ZZ_2$.
\item $(0,1;1,1)$, which gives a fine grading by $\ZZ\times\ZZ_2^2$.
\item   $(0,0;1,1,1,1)$, which gives a fine grading by $ \ZZ_2^4$.
\item $(1,1;\emptyset)$, which gives a fine grading by $\ZZ\times\ZZ_2^3$.
\item $(1,0;1,q_1)$ ($q_1$ as in \eqref{eq:qs}), which gives a fine grading by $\ZZ_2^2\times \ZZ_4$.
\item $(2,0;1)$, which gives a fine grading by $\ZZ_2^5$.
\end{itemize}
\end{itemize}
Let us look at one case in detail. If $m=1$ and $s=0$ then $r=2$ and $\cD=\cQ=M_2(\FF)$. The homogeneous elements in $\cQ$ are $\cup_{i=0}^3\FF q_i$ ($q_0=1$), so that the element $d_1$ can be taken to be $1$ because of the possibility of multiplying by $z$ as in (2.b). The element $d_2$ can be taken to be $q_i$ for some $i$ (multiplying by a nonzero scalar) and, moreover, either $1$ or $q_1$ since there is a graded equivalence of $\cQ$ which sends $q_i$ to $q_1$ if $i\ne 0$. The first possibility: $\Gamma_{(1,0;1,1)}$, provides a grading by $\ZZ_2^4$ which is not fine. Hence we are in the case $\Gamma_{(1,0;1,q_1)}$. The group $\widetilde G$ is generated by  $T$, isomorphic to $\ZZ_2^2$, together with an order two element $g_1$  and another element $g_2$ with square $\deg q_1$ (so that the order of $g_2$ is 4). Hence
$\bar G$ can be identified with $\ZZ_4\times\ZZ_2$ ($\deg q_1=(\bar2,\bar0)$ and $\deg q_2=(\bar0,\bar1)$ generate the subgroup $T$ of $\bar G$) taking   $g_2g_1^{-1}=(\bar1,\bar0)$. A matrix $x\in M_2(\cQ)$ with $q_k$ in the position $(i,j)$ will have degree $(\deg q_k) g_ig_j^{-1}$, so that $M_2(\cQ)$ decomposes as a direct sum of $8$ homogeneous components, all of them of dimension $2$. Here 
\[
-\varphi(x)=\left(\begin{array}{cc}-x_{11}^t&-x_{21}^tq_1\\-q_1x_{12}^t&-q_1x_{22}^tq_1\end{array}\right)
\]
 acts with eigenvalue $-1$ in $M_2(\cQ)_{(\bar0,\bar0)}$ and in $M_2(\cQ)_{(\bar2,\bar0)}$ and splits each   of the remaining homogeneous components of $M_2(\cQ)$ into two pieces, producing a $\ZZ_4\times\ZZ_2^2$-grading on $M_2(\cQ)^-\simeq M_4(\FF)^-$, and hence on $\frsl_4(\FF)=\mathcal{L}$, given by $\mathcal{L}_{e}=0=\mathcal{L}_{(\bar2,\bar0,\bar0)}$, 
\[
 \mathcal{L}_{(\bar2,\bar0,\bar1)}=\left\{
 \left(\begin{smallmatrix}
 a&0&0&0\\0&-a&0&0\\0&0&b&0\\0&0&0&-b
  \end{smallmatrix}\right):a,b\in\FF\right\},
\]
and all the other   13   homogeneous components of dimension 1.
\end{example}

\begin{example} 
An example of an outer fine grading on $\frsl_n(\FF)$ where not only $-\varphi$ is not involutive, but also there are no order two elements in $M'(-\varphi)$, is the $\ZZ_4^3$-fine grading $\Gamma_{(1,0;\,1,q_1,q_2,q_3)}$ on  $\frsl_{8}(\FF)$.
\end{example}

\bigskip

\section{Fine gradings on orthogonal and symplectic Lie algebras}\label{se:BCD}

Recall that the ground field $\FF$ is assumed to be algebraically closed of characteristic zero.

Involutions of the matrix algebra $\cR=M_n(\FF)$ come in two flavors. If $n$ is odd there are only orthogonal involutions, all of them conjugate to the transposition, while if $n$ is even, besides the orthogonal involutions, there appear the symplectic involutions, and all of them are conjugate. If $\varphi$ is an involution of $\cR$, the Lie algebra of skew symmetric elements 
\[
\cK(\cR,\varphi)\bydef\{ x\in \cR : \varphi(x)=-x\}
\]
is isomorphic to the orthogonal Lie algebra $\frso_n(\FF)$ if $\varphi$ is orthogonal, and to the symplectic Lie algebra $\frsp_{2k}(\FF)$ if $n=2k$ and $\varphi$ is symplectic. Moreover, the restriction map
\[
\begin{split}
\Aut(\cR,\varphi)&\longrightarrow \Aut(\cK(\cR,\varphi))\\
\phi\quad&\mapsto \quad \phi\vert_{\cK(\cR,\varphi)},
\end{split}
\]
is a group isomorphism if $n\geq 5$, unless $\varphi$ is orthogonal and $n=6$ or $n=8$ (see \cite[Chapter IX]{Jac79}). If $n=6$, $\frso_6(\FF)$ is isomorphic to $\frsl_4(\FF)$, and for $n=8$, the automorphism group of $\frso_8(\FF)$ contains outer automorphisms of order $3$, due to the phenomenon of triality.

Therefore, with these exceptions, the classification of the fine gradings on $\cK(\cR,\varphi)$ reduces to the classification of fine gradings in $(\cR,\varphi)$. Given such a grading on $(\cR,\varphi)$, we may identify $\cR$ with $\End_\cD(\cV)$, where $\cD=\cQ^{\otimes m}$ for some $m\geq 0$, and $\cV$ is a free right $\cD$-module endowed with a hermitian form $B:\cV\times\cV\rightarrow \cD$. That is, $B$ is sesquilinear, nondegenerate, and  also 
$B(v,w)=\tau\bigl(B(w,v)\bigr)$ for any $v,w\in \cV$, where $\tau$ is the involution on $\cQ^{\otimes m}$ considered in Remark \ref{re:Qm}. Moreover, $\varphi$ is given by the `adjoint' relative to $B$.

Then there is a homogeneous $\cD$-basis $\{v_1,\ldots,v_r,\ldots,v_{r+2s}\}$ of $\cV$ such that the coordinate matrix of $B$ has the following block-diagonal form:
\[ 
M_B= \scriptstyle{
\begin{pmatrix}
d_1&&&&&&&\\
&\ddots&&&&&&\\
&&d_r&&&&&\\
&&&0&1&&&\\
&&&\epsilon&0&&&\\
&&&&&\ddots&&\\
&&&&&&0&1\\
&&&&&&\epsilon&0
\end{pmatrix}}
\]
where $r\geq 0$, $d_1,\ldots,d_r$ are nonzero homogeneous elements in $\cD=\cQ^{\otimes m}$, and either $\epsilon=1$ and $\tau(d_i)=d_i$ for all $i=1,\ldots,r$, if $\varphi$ is orthogonal, or $\epsilon=-1$ and $\tau(d_i)=-d_i$ for all $i=1,\ldots,r$, if $\varphi$ is symplectic. (If $\cD=\FF$ and $\varphi$ is symplectic, this clearly forces $r=0$.)

Moreover, if $\deg (v_i)=g_i$ for $i=1,\ldots,r,\ldots,r+2s$, one has the same relations as in \eqref{eq:gis}. As for type $A$, let  $\widetilde G$ be the abelian group generated by a subgroup $T$ isomorphic to $\ZZ_2^{2m}$ (the grading group of $\cD=\cQ^{\otimes m}$) and elements $g_1,\ldots,g_{r+2s}$, subject only to the relations in \eqref{eq:gis}. The universal group of this fine grading on $(\cR,\varphi)$ is isomorphic to the subgroup $G$ of $\widetilde G$ generated by $T$ and the elements $g_ig_j^{-1}$, $1\leq i,j\leq r+2s$. The free rank of $G$ is exactly $s$ and $G$ is the Cartesian product of a $2$-group that contains $T$ and a free subgroup.

In this situation, attach to this fine grading the sequence $(m,s;d_1,\ldots,d_r)$ as for type $A$ and denote by $\Gamma'_{(m,s;\,d_1,\ldots,d_r)}$ the grading restricted to $\cK(\cR,\varphi)$. 
Then:

\begin{theorem}\label{th:BCD}
Let $\varphi$ be an involution of $\cR=M_n(\FF)$ and assume $n\geq 5$, and $n\ne 6,8$ if $\varphi$ is orthogonal. Then any fine grading on $\cK(\cR,\varphi)$ is equivalent to a grading $\Gamma'_{(m,s;\,d_1,\ldots,d_r)}$ as above for some $m,r,s$ such that $n=2^m(r+2s)$, and homogeneous elements $d_1,\ldots,d_r\in\cQ^{\otimes m}$, except for $s=0$, $r=2$, and $\FF d_1=\FF d_2$.

Two such gradings $\Gamma'_{(m,s;\,d_1,\ldots,d_r)}$ and $\Gamma'_{(m',s';\,d'_1,\ldots,d'_{r'})}$ are equivalent if and only if $m=m'$, $s=s'$, $r=r'$ and there  {exist} a graded equivalence $\Phi\colon\cQ^{\otimes m}\rightarrow \cQ^{\otimes m}$, a permutation $\sigma\in      {S_r}$ and a nonzero homogeneous element $z\in\cQ^{\otimes m}$ such that $\tau(z)=z$ and $\Phi(zd_i)\in \FF d'_{\sigma(i)}$ for all $i=1,\ldots,r$. 
\end{theorem}

For details see \cite{Eld10} or \cite[Chapter 3]{EK13}. The assumption $\tau(z)=z$ does not appear in \cite[Theorem 5.2]{Eld10} because the involution $\tau$ on $\cQ^{\otimes m}$ is not fixed there.

\begin{example}[Fine gradings on $\frso_5(\FF)$]

In this case, if $5=2^m(r+2s)$ then $m=0$ and $r=1,3,5$. In particular $\cD=\FF$ and there are three different  nonequivalent fine gradings:
\begin{itemize}
\item $r=1$, which gives the Cartan grading of $\frso_5(\FF)$ by $\ZZ^2$.
\item $r=3$, which gives a fine grading by $\ZZ\times \ZZ_2^2$.
\item $r=5$, which gives a fine grading by $\ZZ_2^4$.
\end{itemize}
\end{example}

In general, for $\frso_{2k+1}(\FF)$ (type $B$), there are exactly ${k}+1$ {nonequivalent} fine gradings,    
{ since the only possibilities are $M_B=\diag{\{I_r,I_2,\dots,I_2\}}$ and the grading is determined by the number $s$ of $I_2$-blocks ($s\in\{0,\dots,k\}$), with universal grading group $\ZZ^s\times \ZZ_2^{2(k-s)}$}.

\begin{example}[Fine gradings on $\frsp_6(\FF)$]

In this case $6=2^m(r+2s)$, so the following nonequivalent possibilities appear (with $q_3$ as in \eqref{eq:qs}):
\begin{itemize}
\item $(0,3;\emptyset)$, which gives the Cartan grading by $\ZZ^3$.

\item $(1,1;q_3)$, which gives a fine grading by $\ZZ\times\ZZ_2^2$.

\item $(1,0;q_3,q_3,q_3)$, which gives a fine grading by $\ZZ_2^4$.
\end{itemize}
\end{example}

\begin{example}[Fine gradings on $\frsp_8(\FF)$]\label{ex_c4}

In this case $8=2^m(r+2s)$, so there are 7 nonequivalent possibilities:
\begin{itemize}
\item $(0,4;\emptyset)$, which gives the Cartan grading by $\ZZ^4$.

\item $(1,2;\emptyset)$, which gives a fine grading by $\ZZ^2\times\ZZ_2^2$.
\item $(1,1;q_3,q_3)$, which gives a fine grading by $\ZZ\times\ZZ_2^3$.
\item $(1,0;q_3,q_3,q_3,q_3)$, which gives a fine grading by $\ZZ_2^5$.
\item $(2,1;\emptyset)$, which gives a fine grading by $\ZZ\times\ZZ_2^4$.
\item $(2,0;1\otimes q_3,q_3\otimes 1)$, which gives a fine grading by $\ZZ_4\times\ZZ_2^3$.
\item $(3,0;\emptyset)$, which gives a fine grading by $ \ZZ_2^6$.

\end{itemize}
\end{example}

\smallskip

\begin{remark}\label{re:D4}
The situation for $\frso_8(\FF)$ is more complicated. If $\varphi$ is an orthogonal involution of $\cR=M_8(\FF)$, and we identify $\Aut(\cR,\varphi)$ with a subgroup of $\Aut(\frso_8(\FF))$ (by restriction), then $\Aut(\cR,\varphi)$ has index three in $\Aut(\frso_8(\FF))$. It turns out that {any} maximal quasitorus of $\Aut(\frso_8(\FF))$ satisfies one of the following possibilities (see \cite[Theorem 6.7]{Eld10}):
\begin{itemize}
\item Either it is conjugate to a maximal quasitorus of $\Aut(\cR,\varphi)$. There are $15$ such possibilities up to conjugation in $\Aut(\cR,\varphi)$, but two of them are conjugate in $\Aut(\frso_8(\FF))$, so we obtain here $14$ fine gradings up to equivalence.

\item Or it contains an outer automorphism $\theta$ of order $3$. There are, up to conjugation, only two such automorphisms. The centralizer of $\theta$ in $\Aut(\frso_8(\FF))$ is $\langle\theta\rangle\times H$, where $H$ is, up to isomorphism, the simple group of type $G_2$ in one case and $\PGL_3(\FF)$ in the other case. The maximal quasitori of $G_2$ are well-known, while the maximal quasitori of $\PGL_3(\FF)$ correspond to the inner fine gradings on $\frsl_3(\FF)$ and there are only two of them, according to Theorem \ref{th:A}: the Cartan grading and the Pauli grading. As a consequence, there are three more nonequivalent fine gradings on $\frso_8(\FF)$ with universal groups $\ZZ^2\times\ZZ_3$, $\ZZ_2^3\times\ZZ_3$ and $\ZZ_3^3$. 
\end{itemize}
\end{remark}

\bigskip


\section{Fine gradings on the exceptional simple Lie algebras}\label{se:exceptional}

As in the previous two sections, the ground field $\FF$ is assumed here to be algebraically closed of characteristic zero.

The simple Lie algebra $\cL$ of type $G_2$ (respectively $F_4$) is, up to isomorphism, the Lie algebra of derivations of the algebra of octonions $\OO$ (respectively, of the Albert algebra $\Albert$, i.e., the simple exceptional  Jordan algebra of the hermitian matrices of order 3 with coefficients in $\OO$), and any automorphism of $\OO$ (resp., $\Albert$) induces, by conjugation, an automorphism of $\cL$, thus giving an isomorphism of the automorphism groups. Thus the problem of classifying fine gradings, up to equivalence, on $\cL$ reduces to the same problem on the smaller algebra $\OO$ (resp., $\Albert$). (See \cite[Chapter 5]{EK13} and references {therein}.)

\begin{theorem}\label{th:GF} \null\quad
\begin{itemize}
\item Up to equivalence, there are two fine gradings on the simple Lie algebra of type $G_2$: the Cartan grading by $\ZZ^2$ and a grading by $\ZZ_2^3$ in which $\cL_{(\bar 0,\bar 0,\bar 0)}=0$ and $\cL_\alpha$ is a Cartan subalgebra of $\cL$ for any $0\ne \alpha\in\ZZ_2^3$. This grading is induced by the natural $\ZZ_2^3$-grading on $\OO$ obtained by constructing $\OO$ from the ground field in three steps by means of the Cayley-Dickson doubling process.

\item Up to equivalence, there are four fine gradings on the simple Lie algebra of type $F_4$: the Cartan grading by $\ZZ^4$ and gradings by $\ZZ\times\ZZ_2^3$, $\ZZ_2^5$ and $\ZZ_3^3$. For the last grading,  $\cL_\alpha\oplus\cL_{-\alpha}$ is a Cartan subalgebra of $\cL$ for any $0\ne \alpha\in\ZZ_3^3$.
\end{itemize}
\end{theorem}

\smallskip

For the simple Lie algebra of type $E_6$, the fine gradings have been classified in \cite{DV_E6}. Up to equivalence, there are $14$ different such gradings.  If the grading is produced by a maximal quasitorus $M$ of the group of inner automorphisms and $M$ is not a maximal torus, then it contains either an elementary 2-group of type $\ZZ_2^3$ or an elementary 3-group of type $\ZZ_3^2$ (two possibilities here). The knowledge of the three centralizers allows to obtain the possible maximal quasitori, living inside either $\ZZ_2^3\times \PGL(3)$,  $\ZZ_3^2\times \PGL(3)$, or $\ZZ_3^2\times G_2$. This gives 4 maximal quasitori, producing fine gradings by the universal groups $\ZZ_2^3\times \ZZ^2$, $\ZZ_2^3\times \ZZ_3^2$,    $\ZZ^2\times \ZZ_3^2$, and $ \ZZ_3^4$. Otherwise $M$ contains  outer automorphisms.
If $M$ contains an order two outer automorphism, this automorphism fixes a subalgebra of type either $C_4$ or $F_4$ and the grading comes from extending either a fine grading on $C_4$ (seven possibilities here, according to Example~\ref{ex_c4}) or a fine grading on $F_4$ (4 possibilities by Theorem~\ref{th:GF}, three of them containing also automorphisms fixing $C_4$). On the other hand, if $M$ contains outer automorphisms but none of them has order two, then the quasitorus $M$ is isomorphic to $\ZZ_4^3$: an outer automorphism fixes a subalgebra isomorphic to $\frsl_4(\FF)\oplus \frsl_2(\FF)$ and the restriction of the fine grading to $\frsl_4(\FF)$ is just the Pauli $\ZZ_4^2$-grading.  

\smallskip

The classification for $E_7$ and $E_8$ can be derived from recent work of Jun Yu \cite{Yu14}. Yu classifies the conjugacy classes of the closed abelian subgroups $F$ of the compact real simple Lie groups $G$ satisfying the condition $\dim\frg_0^F=\dim F$, where $\frg_0$ is the Lie algebra of $G$ and $\frg_0^F$ is the subalgebra of fixed elements by the action of $F$. This class of groups presents some nice functorial properties. In particular, the maximal finite abelian subgroups are among these subgroups. The close relationship between compact Lie groups and complex reductive algebraic groups allows, in principle, to extract from \cite{Yu14} the list of the conjugacy classes of the maximal quasitori of the automorphism groups of the simple exceptional Lie algebras over $\CC$. This gives the classification of the equivalence classes of fine gradings in these algebras. The results over $\CC$ can then be transferred to any algebraically closed field of characteristic zero \cite{Eldpr}.

{Then} it turns out that the tentative list in \cite[Figure 6.2]{EK13} is complete. 
 Up to equivalence, there are again $14$ fine gradings both on $E_7$ and on $E_8$, although only some of them form natural families in $E_6$, $E_7$ and $E_8$.
 This list contains the  universal grading groups of the fine gradings on the simple Lie algebras of types $E_6$, $E_7$ and $E_8$  together with    a convenient model in each case which stresses how some of these gradings appear in natural families.  These models of the fine gradings on the simple Lie algebras of type $E$ are thoroughly discussed in \cite{DE_Note}. This work emphasizes the role of the nonassociative algebras in the gradings, describing them  by using not only the famous unified construction of the exceptional Lie algebras by Tits, but also   constructions based on symmetric composition algebras (specially relevant for explaining the $\ZZ_3$-gradings) and the Kantor and Steinberg contructions of Lie algebras out of structurable algebras (related to $\ZZ$-gradings with more than three pieces).
 
\smallskip

The fine gradings by finite groups on the simple Lie algebra of type $E_8$ have been independently   classified in \cite{MADsE8}. This problem is equivalent to the computation of the conjugacy classes of the maximal abelian finite subgroups of the simple algebraic group $E_8$.  These maximal quasitori  (and hence the universal grading groups of the related fine gradings) are isomorphic to either $\ZZ_3^5$, $\ZZ_6^3$, $\ZZ_2^9$, $\ZZ_2^8$, $\ZZ_4^3\times\ZZ_2^2$, $\ZZ_4\times\ZZ_2^6$, or $\ZZ_5^3$. One of the main tools used for this classification is the Brauer invariant of the irreducible modules for graded semisimple Lie algebras introduced in \cite{EKgr}. This paper studies conditions for a module to be graded in a way compatible with a given  grading on the Lie algebra. 
The approach is thus quite different to the one in \cite{Yu14}. 
These fine gradings on simple Lie algebras by finite groups are remarkable because their behaviour is completely different to the one of the root space decomposition (they may be considered just at the other end of the spectrum of fine gradings): for instance, every nonzero homogeneous element is (ad-)semisimple, which allows to choose bases formed by semisimple elements. (In general, gradings are closely related to the problem of a suitable choice of basis. Recall the relationship Chevalley basis--Cartan grading.) 

 A remarkable grading in the above list is the $\ZZ_5^3$-grading on $E_8$, because such $5$-symmetry is a particular fact of $E_8$. Besides, it is one of the so-called
  \emph{Jordan gradings}  \cite{Eld_Jordan},
 as well as the above fine $\ZZ_2^3$-grading on the simple Lie algebra of type $G_2$ and the fine $\ZZ_3^3$-grading on the simple Lie algebra of type $F_4$. These three gradings satisfy that  every nonzero homogeneous component has dimension $2$ and is contained in a Cartan subalgebra. 
 The fine $\ZZ_5^3$-grading on $E_8$ has not gone unnoticed. The interested reader can consult the notes of Kostant's talk \cite{Kost}, which deals with this and other gradings, like Dempwolf's decomposition of the Lie algebra of type $E_8$ as a sum of $31$ pieces, all of them Cartan subalgebras, which is a $\ZZ_2^5$-grading, obtained as a  coarsening of the fine $\ZZ_2^8$-grading.  
 
\smallskip

Note that the problem of the classification of fine gradings by finite groups is a key piece of the puzzle of the (general) classification,  because if the universal group   is infinite, then the grading on the Lie algebra induces a grading by a not necessarily reduced root system \cite{Eld_Fine} and it is determined by a fine grading on the coordinate algebra of the grading by the root system. Associative, alternative, Jordan or structurable algebras appear as coordinate algebras. In a sense, the classification of the fine gradings whose associated quasitori are not finite is reduced to the classification of some fine gradings on certain nonassociative algebras.

\bigskip


\section{Modular case}\label{se:modular}

In this section, the ground field $\FF$ is assumed to be algebraically closed of characteristic not two.

The gradings on a finite dimensional algebra $\cA$ are no longer given by means of common eigenspaces for the action of an (abelian) diagonalizable subgroup of $\Aut(\cA)$. A different approach is needed.

Given a grading $\Gamma$ as in \eqref{eq:grading} by the group $G$, consider the map:
\[
\begin{split}
\eta\colon\cA&\longrightarrow \cA\otimes \FF G\\
x_g&\mapsto\ x_g\otimes g
\end{split}
\]
for any $g\in G$ and $x_g\in\cA_g$, where $\FF G$ denotes the group algebra of $G$. Then $\eta$ is both a homomorphism of algebras and a map that provides $\cA$ with the structure of a comodule for the Hopf algebra $\FF G$. The map $\eta$ is then called a \emph{comodule algebra map}.

Conversely, given such a map $\eta$, $\cA$ is graded by $G$ with
\[
\cA_g=\{ x\in \cA: \eta(x)=x\otimes g\}
\]
for any $g\in G$. In a way, this means that $\cA_g$ is the eigenspace for $\eta$ with eigenvalue $g$. Thus, gradings by $G$ on $\cA$ correspond bijectively with the comodule algebra maps $\cA\rightarrow \cA\otimes \FF G$.

But any comodule algebra map $\eta$ induces a \emph{generic} automorphism of algebras over $\FF G$:
\begin{equation}\label{eq:AFG}
\begin{split}
\cA\otimes \FF G&\longrightarrow \cA\otimes \FF G\\
x\otimes h\ &\mapsto\ \eta(x)  h 
\end{split}
\end{equation}
so that $x_g\otimes h\mapsto x_g\otimes gh$ for any $g\in G$ and $x_g\in\cA_g$. All the information on $\Gamma$ is contained in this single automorphism.

More generally, a comodule algebra map $\eta:\cA\rightarrow \cA\otimes \FF G$ defines a homomorphism of affine group schemes:
\[
\rho: G^D\longrightarrow \AAut(\cA),
\]
where the `Cartier dual' $G^D$ is the affine group scheme (i.e., the representable functor from the category of unital asociative commutative algebras over $\FF$ into the category of groups) such that \[
G^D(R)=\Hom_{\textrm{alg}}(\FF G,R)\simeq \Hom_{\textrm{groups}}(G,R^\times),
\] 
and $\AAut(\cA)$ is the affine group scheme whose $R$-points are the automorphisms of the $R$-algebra $\cA\otimes R$: $\AAut(\cA)(R)=\Aut_{\textrm{$R$-alg}}(\cA\otimes R)$. The behavior  of $\rho$ on homomorphisms is the natural one. For each unital associative commutative $\FF$-algebra $R$, the map $\rho_R$ is defined as follows:
\[
\begin{split}
\rho_R: G^D(R)=\Hom_{\textrm{alg}}(\FF G,R)&\longrightarrow \AAut(\cA)(R)\\
f:\FF G\rightarrow R\ &\mapsto\quad \rho_R(f):\cA\otimes R\rightarrow \cA\otimes R\\
                       &\hspace{60pt}  x_g\otimes r\mapsto x_g\otimes f(g)r.
\end{split}
\]
Conversely, if $\rho:G^D\rightarrow \AAut(\cA)$ is a homomorphism of affine group schemes (i.e., a natural transformation), $\rho_{\FF G}(\id)$ is an automorphism of $\FF G$-algebras $\cA\otimes \FF G\rightarrow \cA\otimes\FF G$ as in \eqref{eq:AFG}, which induces a comodule algebra map by composition:
\[
\cA\simeq \cA\otimes 1\hookrightarrow \cA\otimes \FF G\xrightarrow{\rho_{\FF G}(\id)}\cA\otimes\FF G.
\]

The conclusion is that gradings by $G$ on $\cA$ correspond bijectively to homomorphisms of affine group schemes $G^D\rightarrow \AAut(\cA)$.

In other words, to work in prime characteristic, we have to substitute the group of characters $\hat G$ by the Cartier dual $G^D$, and the algebraic group $\Aut(\cA)$ by the affine group scheme $\AAut(\cA)$. With this in mind, Theorem \ref{th:gradings_MADs} remains valid: the classification of the fine gradings up to equivalence corresponds to the classification of the maximal quasitori in $\AAut(\cA)$ up to conjugation by elements in $\Aut(\cA)$. (See \cite{EK13} for details.)

In particular, if two algebras have isomorphic affine group schemes of automorphisms, we can transfer the problem of classification from one algebra to the other. The outcome is that Theorems \ref{th:A}, \ref{th:BCD}, and \ref{th:GF} remain valid in the modular case if we change $\frsl_n(\FF)$ by $\frpsl_n(\FF)$  ($=[\cR,\cR]/(Z(\cR)\cap [\cR,\cR])$ for $\cR=M_n(\FF)$) with a couple of exceptions:
\begin{itemize}
\item The $\ZZ_3^3$-grading on the simple Lie algebra of type $F_4$ does not exist in characteristic $3$.

\item Also in characteristic $3$, the automorphism group scheme $\AAut(\frpsl_3(\FF))$ is not isomorphic to the group scheme of automorphisms and antiautomorphisms of $M_3(\FF)$, but to the group scheme of automorphisms of the octonions! Hence in this case there are only two fine gradings on $\frpsl_3(\FF)$ with universal groups $\ZZ^2$ and $\ZZ_2^3$. Moreover, in this situation there is no simple Lie algebra of type $G_2$. (See \cite{CE16} for some related results.)
\end{itemize}

For the simple Lie algebra $\cL$ of type $D_4$ (see \cite{EK_D4}), in characteristic $3$ all the fine gradings are obtained by restriction of fine gradings in $(M_8(\FF),\text{t})$, where $\text{t}$ denotes the transpose involution (so there are $14$ fine gradings up to equivalence), while if the characteristic is $>3$, the results in characteristic $0$ remain valid, but with a different proof  that relies in the fact that $\AAut(\cL)$ is isomorphic to the affine group scheme of automorphisms of   certain algebraic structure called \emph{trialitarian algebra}.
The general philosophy is to find a simpler object sharing the  affine group scheme of automorphisms with the Lie algebra under study.

\smallskip

For information on gradings on some simple modular Lie algebras of Cartan type, the reader may consult \cite[Chapter 7]{EK13}.

\smallskip

We finish this survey with the following

\medskip

\noindent\textbf{Open problem:} Classify the fine gradings, up to equivalence, on the simple Lie algebras of types $E_6$, $E_7$ and $E_8$ over fields of prime characteristic ($\ne 2$).

\bigskip


\end{document}